\documentclass[]{article}
\usepackage{geometry}
\usepackage{graphicx}
\usepackage{amssymb}
\usepackage{amsmath}
\usepackage{amsthm}
\usepackage{subfigure}
\usepackage{wrapfig}
\usepackage{stmaryrd}
\usepackage{mathrsfs}
\usepackage{bbm}
\usepackage{thm-restate}
\usepackage{quiver}
\usepackage{hyperref}
\title{On the profinite distinguishability of hyperbolic Dehn fillings of finite-volume 3-manifolds}
\author{Paul Rapoport} 
\date{\today}

\begin{document}
\maketitle

\newcommand{\fr}[1]{\mathfrak{#1}}
\newcommand{\bb}[1]{\mathbb{#1}}
\newcommand{\scr}[1]{\mathscr{#1}}
\newcommand{\ds}{\displaystyle}

\newtheorem{theorem}{Theorem}[section]
\newtheorem*{theorem*}{Theorem}
\newtheorem{corollary}{Corollary}[theorem]
\newtheorem{lemma}[theorem]{Lemma}
\newtheorem{lemma*}{Lemma}
\newtheorem*{remark}{Remark}
\newtheorem*{conjecture*}{Conjecture}
\newtheorem*{question*}{Question}
\newtheorem{thmA}{Theorem}
\renewcommand{\thethmA}{\Alph{thmA}}
\newtheorem{corA}[thmA]{Corollary}
\newtheorem{definition}[theorem]{Definition}
\newtheorem{proposition}[theorem]{Proposition}
\newcommand{\forwards}{$(\Rightarrow)$}
\newcommand{\reverse}{$(\Leftarrow)$}
\newcommand{\bicond}{$\Leftrightarrow$}

\newcommand{\speq}{\supseteq}
\newcommand{\sbeq}{\subseteq}
\newcommand{\spneq}{\supsetneq}
\newcommand{\sbneq}{\subsetneq}

\newcommand{\half}{\frac{1}{2}}
\newcommand{\invn}[1]{\frac{1}{#1}}

\newcommand{\nonz}{\backslash \{0\}}
\newcommand{\nonid}{\backslash \{1\}}
\newcommand{\eps}{\epsilon}

\newcommand{\inv}{^{-1}}

\newcommand{\cd}{\cdot}
\newcommand{\cds}{\cdots}

\newcommand{\nat}{\mathbb{N}}
\newcommand{\z}{\mathbb{Z}}
\newcommand{\q}{\mathbb{Q}}
\newcommand{\real}{\mathbb{R}}
\newcommand{\complex}{\mathbb{C}}
\newcommand{\finfield}[1]{\mathbb{F}_{#1}}
\newcommand{\finfieldcls}[1]{\overline{\mathbb{F}}_{#1}}
\newcommand{\nneg}{_{\geq 0}}

\newcommand{\cpx}{\complex}

\newcommand{\map}{\rightarrow}
\newcommand{\inj}{\hookrightarrow}
\newcommand{\sur}{\twoheadrightarrow}
\newcommand{\iso}{\stackrel{\sim}{\rightarrow}}

\newcommand{\adjoin}[1]{\left[#1\right]}
\newcommand{\labelmap}[1]{\stackrel{#1}{\rightarrow}}
\newcommand{\labelinj}[1]{\stackrel{#1}{\hookrightarrow}}
\newcommand{\labelsur}[1]{\stackrel{#1}{\twoheadrightarrow}}

\newcommand{\shortexact}[3]{0 \rightarrow #1 \hookrightarrow #2 \twoheadrightarrow #3 \rightarrow 0}
\newcommand{\labelshortexact}[5]{0 \rightarrow #1 \stackrel{#2}{\hookrightarrow} #3 \stackrel{#4}{\twoheadrightarrow} #5 \rightarrow 0}

\newcommand{\shortexactone}[3]{1 \rightarrow #1 \hookrightarrow #2 \twoheadrightarrow #3 \rightarrow 1}
\newcommand{\labelshortexactone}[5]{1 \rightarrow #1 \stackrel{#2}{\hookrightarrow} #3 \stackrel{#4}{\twoheadrightarrow} #5 \rightarrow 1}

\newcommand{\gen}[1]{\langle #1 \rangle}
\newcommand{\ncls}[1]{\gen{\gen{#1}}}

\newcommand{\nsgp}{\lhd}
\newcommand{\nsgpeq}{\unhld}

\newcommand{\imm}{\looparrowright}
\newcommand{\bdd}{\partial}

\newcommand{\action}{\curvearrowright}
\newcommand{\itvl}{[0, 1]}
\newcommand{\oball}[2]{B_{#1}(#2)}
\newcommand{\cball}[2]{\bar{B}_{#1}(#2)}

\newcommand{\infser}[3]{\{#1_#2\}_{#2 = #3}^\infty} 
\newcommand{\infsum}[2]{\displaystyle\sum_{#1 = #2}^\infty}
\newcommand{\infprod}[2]{\displaystyle\prod_{#1 = #2}^\infty}
\newcommand{\infcup}[2]{\displaystyle\bigcup_{#1 = #2}^\infty}
\newcommand{\infcap}[2]{\displaystyle\bigcap_{#1 = #2}^\infty}

\newcommand{\pderiv}[2]{\frac{\partial #1}{\partial #2}}
\section{Introduction}
\label{intro}
One of the foundational lines of research in modern geometry has been the study of $3$-manifold topology; some relevant recent progress has been made primarily through teasing out differences between different manifolds through appeal to their fundamental groups. For some examples, see \cite{wise2009research}, \cite{agol2013virtual}, \cite{wilton2017distinguishing}, and \cite{przytycki2018mixed}.\\
It is well known that fundamental groups of finite-volume hyperbolic $3$-manifolds are finitely presented, and that the fundamental groups of finite-volume hyperbolic $3$-manifolds are residually finite (see \cite{malcev1940isomorphic}). Thus it seems like a natural line of inquiry to try to distinguish the fundamental groups of finite-volume $3$-manifolds by looking at their finite quotients, since finite quotients of fundamental groups correspond to finite-sheeted regular coverings, so that commensurability of $3$-manifolds coincides with commensurability of their fundamental groups.\\
Tying together all of the information about finite quotients of fundamental groups through an inverse limit leads us to the concept of the profinite completion of the group (see Definition \ref{profin-def}), and to profinite rigidity (see Definition \ref{profin-rigid}); this thus follows in the path laid down by \cite{noskov1979infinite}, \cite{reid2013profinite}, and \cite{wilton2019profinite}, among many others. In particular, Wilton and Zaleskii in \cite[Thm 8.4]{wilton2017distinguishing} show that the profinite completion of a geometric $3$-manifold group determines its geometry, and further, in \cite{wilton2019profinite}, they show that the profinite completion of a $3$-manifold group also determines the JSJ decomposition of the manifold.\\
There are two natural notions of profinite rigidity: relative and absolute. Relative profinite rigidity deals with whether within some family of groups, a given pair of distinct groups can be distinguished by their profinite completions, while absolute profinite rigidity is more general, dealing with whether a given group is profinitely distinguishable from every other residually finite group. (For more precise definitions, see Definition \ref{profin-rigid} later on.) Since it is a more powerful and general property, we might expect that there have been some attempts to construct a few $3$-manifold groups which are provably profinitely rigid in the absolute sense. In Bridson, McReynolds, Reid, and Spitler's breakthrough paper \cite{bridson2018absolute}, the authors do just that, producing the only known examples of absolutely profinitely rigid hyperbolic $3$-manifold groups, through careful examination of certain arithmetic lattices in $\mathrm{PSL}(2, \cpx) \cong \mathrm{Isom}^+(\bb{H}^3)$ and the employment of strong assumptions, though thus far there are only finitely many known absolutely profinitely rigid $3$-manifold groups, and very little else is known about the absolute case. Accordingly, I have restricted my work to the realm of relative profinite rigidity within the category of $3$-manifold groups, since its study has more work done for me to build on and more tools available for me to use. The powerful work of Agol in \cite{agol2013virtual} and Wise in \cite{wise2009research} in particular provided many tools for the study of finite covers of $3$-manifolds, opening up in turn the ability to make more progress on the study of profinite rigidity within the realm of $3$-manifold groups, and my work follows this side of the branch: my use of the non-Haken assumption here comes from the foundational work of \cite{culler1983varieties}.\\
The study of finite-volume hyperbolic $3$-manifolds is a central area within $3$-manifold topology: among $3$-manifolds, the possession of a hyperbolic structure is the ``generic'' case, as we see by Thurston's Hyperbolic Dehn Surgery Theorem in \cite[Theorem 2.6]{thurston1982three} and which is explored in more formal and precise detail by Maher in \cite{maher2010random}. We recall further that the study of Dehn fillings is fundamental to $3$-manifold topology, and that the volumes of the Dehn fillings $M_{p/q}$ - notation which is defined later on in Definition \ref{dehn-filling-defn} - of a hyperbolic $3$-manifold $M$ are strictly less than and converge to $vol(M)$,\footnote{In fact, the order type of the set of all volumes of Dehn fillings of all $3$-manifolds is $\omega^\omega$!} by \cite[Thurston's Theorem]{gromov1981hyperbolic}. Additionally, Mostow's Rigidity Theorem, here Theorem \ref{mostow-rigidity}, says that geometric invariants of a complete and finite-volume hyperbolic $3$-manifold $M$, such as volume, are topological invariants, and thus give invariants of their fundamental groups. Accordingly, we have even further reason to suspect that the study of finite-volume hyperbolic $3$-manifolds might be fruitful. Mostow rigidity implies that the fundamental groups $\pi_1(M_{p/q})$ are different, and it becomes a natural question whether the profinite completions of these groups are also distinguishable. This is a question that I answer in the affirmative in cofinitely many cases.\\
It is an open question as to whether all finite-volume hyperbolic $3$-manifold groups are absolutely profinitely rigid, and furthermore, it is also currently unknown whether there exists some pair of non-isometric hyperbolic manifolds that are not profinitely distinguishable. However, my work presented here has helped to further explore this line of inquiry by ruling out a class of possibilities. In particular, we have the following theorem, proved in Section \ref{sec-main}, with the precise definition of the notation $|\chi^I_\cpx(\Gamma)|$ I use for character varieties in the theorem below found in Definition \ref{charvar-notation}:

\begin{restatable}{thmA}{mainone}
\label{mainone}
Let $\Gamma$ be any finitely generated residually finite group with $|\chi^I_\cpx(\Gamma)| < \infty$, and let $M$ be an oriented finite-volume hyperbolic 3-manifold with a single cusp. Then for hyperbolic Dehn fillings with all but finitely many choices of surgery coefficient $M_{m/n}$ with $\Lambda = \pi_1(M_{m/n})$, we get that $\widehat{\Gamma} \not\cong \widehat{\Lambda}$.
\end{restatable}

Our journey takes us through representation theory, as well: it turns out to be easier to look at the $\mathrm{SL}(2, k)$-representations of $3$-manifold groups over careful choices of field $k$ rather than at the groups themselves, and in turn at the character variety of a given representation rather than at the representation itself. By a theorem of Culler and Shalen, \cite[Proposition 1.5.2]{culler1983varieties}, which here is Theorem \ref{D}, a point in the character variety of a $3$-manifold group picks out an irreducible representation up to conjugacy, and by a result of \cite{culler1983varieties} - Theorem \ref{A} here - the character variety of a non-Haken $3$-manifold group is in fact finite. I owe a debt of gratitude to Bridson, McReynolds, Reid, and Spitler \cite{bridson2018absolute}, who have used representation theory to think about profinite distinguishability of specific $3$-manifold groups in a different way.\\
The ``special sauce'' here is my use of model theory as described by Marker \cite{marker1996introduction}, employing a Lefschetzian transfer principle rather than algebraic geometry. This framework ensures that I can pass back and forth between $\mathrm{SL}(2, \cpx)$ and $\mathrm{SL}(2, \finfieldcls{p})$-representations as needed, which is crucial, as both have desirable properties: $\mathrm{SL}(2, \finfieldcls{p})$ is locally finite, which we find useful in controlling profinite extensions of maps in both Corollary \ref{local-finiteness} and Lemma \ref{profin-same-size-p}, but by contrast, representations into $\mathrm{SL}(2, \cpx)$ are both much better understood and more directly connected to more concrete geometric applications. In particular, any finite-volume hyperbolic manifold corresponds naturally to a finite-covolume lattice within $\mathrm{(P)SL}(2, \cpx)$, and this gives us a canonical representation.\\
While this thesis was in preparation, \cite{liu2020finitevolume} appeared on the ArXiV, in which Liu proves a more general version of Corollary \ref{cmainone} using sophisticated methods more closely hewing to orthodox geometric group theory, and which as a consequence concerns itself purely with $3$-manifold groups; by contrast, my borrowing from model theory has the major advantages of being a more elementary approach while permitting substantively more generality of hypothesis.

\section{Preliminaries}
\label{setup}

\subsection{Hyperbolic Geometry and Dehn Surgery}
For a very readable treatment of foundational concepts in $3$-manifold topology, see Hatcher \cite{hatcher-notes}.

\begin{definition}
Let $M$ be a $3$-manifold, $S$ a compact surface properly embedded in $M$. Suppose there exists some disk $D \sbeq M$ such that $D \cap S = \bdd D$, with the intersection being transverse. If $\bdd D$ bounds no disk in $S$, then we call $D$ a \emph{nontrivial compressing disk}, and $S$ is a \emph{compressible surface}. Otherwise, if $S$ neither has a nontrivial compressing disk nor is an embedded $2$-sphere, then $S$ is an \emph{incompressible surface}.
\end{definition}

\begin{definition}
Let $M$ be a connected $3$-manifold. We call $M$ \emph{prime} if there exist no $3$-manifolds $N_1, N_2 \not \cong S^3$ such that $N_1 \# N_2 \cong M$. We call $M$ \emph{irreducible} if every $S^2 \sbeq M$ bounds a $3$-ball.
\end{definition}

\begin{lemma}
\cite[Lemma 1]{milnor1962unique} Let $M$ be a prime $3$-manifold. Then either $M$ is irreducible or $M \cong S^2 \times S^1$.
\end{lemma}

\begin{definition}
Let $M$ be a compact, orientable, prime $3$-manifold. Then we say that $M$ is \emph{Haken} if it has at least one properly embedded two-sided incompressible surface; if it has none, we call it \emph{non-Haken}.
\end{definition}

Synonymously, and perhaps more suggestively, Haken manifolds are also called \textit{sufficiently large}, because the surfaces from the definition are seams along which we can repeatedly cut the manifold, getting simpler pieces each time - which are still Haken \cite[Theorem 13.3]{hempel2004three} - until at last we have a disjoint union of $3$-balls. For example, consider the three-torus $S^1 \times S^1 \times S^1$, which we may think of as a solid cube $[0, 1] \times [0, 1] \times [0, 1]$ with opposite faces identified. Cutting along the torus $[0, 1] \times [0, 1] \times \{\invn{2}\}$ yields a torus cross an interval, which we may think of as a slightly thickened $S^1 \times S^1$, and then cutting along an incompressible annulus yields a solid torus $S^1 \times B^2$. Finally, cutting along a disk yields a $3$-ball.\\
This then permits us to prove complex statements about an entire class of $3$ manifolds by induction on the length of the Haken hierarchy: we show that a property holds of the $3$-ball, and then that the property is preserved through gluing along an incompressible surface.

\begin{definition}
\label{dehn-filling-defn}
Let $M$ be a $3$-manifold such that $\bdd M$ consists of a single torus $T \cong S^1 \times S^1$, with $H_1(\bdd T)$ generated by choices of longitude $l$ and meridian $m$. For $\frac{p}{q} \in \q \cup \infty$, the \emph{Dehn filling of $M$ along $T$ with slope $\frac{p}{q}$} is given by $M \cup T'_{p/q}$, for $T'$ a solid torus $T' \cong B^2 \times S^1$, where the union is a gluing along $T$ such that the meridian of $T'$ maps to a corresponding curve in $\bdd T$ homotopic to $q\cdot[m] + p\cdot[l]$. We denote the resulting manifold by $M_{p/q}$.
\end{definition}

\begin{definition}
Let $M$ be an orientable finite-volume hyperbolic $3$-manifold, and let $T_1, \cdots, T_n$ be its boundary components on truncating each cusp. Then for $\frac{p_1}{q_1}, \frac{p_2}{q_2}, \cdots, \frac{p_n}{q_n} \in \q \cup \infty$, the \emph{hyperbolic Dehn surgery with slopes $(\frac{p_1}{q_1}, \frac{p_2}{q_2}, \cdots, \frac{p_n}{q_n})$} is given by gluing in $n$ solid tori, one for each component of the boundary, where gluing are given by mapping the meridian of each $T'_i$ to the curve $q_i[l_i] + p_i[m_i]$ of $T_i$.
\end{definition}

\begin{theorem}
(Thurston's hyperbolic Dehn surgery theorem) \cite[Theorem 2.6]{thurston1982three} Let $M$ be a hyperbolic $3$-manifold with $n$ cusps, and let $M(p_1/q_1, p_2/q_2, \cdots p_n/q_n)$ be the manifold obtained through applying a hyperbolic Dehn filling to each of the $n$ cusps with slopes $p_i/q_i$. Then if each of the $p_i/q_i$ differs from finitely many exceptional slopes, $M(p_1/q_1, p_2/q_2, \cdots p_n/q_n)$ is also hyperbolic.
\end{theorem}

\begin{definition}
Let $M$ be a connected Riemannian manifold. We say that $M$ is a \emph{hyperbolic $3$-manifold} if it is complete and everywhere locally isometric to the hyperbolic $3$-space $\bb{H}^3$.
\end{definition}

\begin{remark}
Many knot complements $S^3 \backslash K$ have hyperbolic structure; a notable family of exceptions is the set of torus knots. For example, the figure-eight knot has a complement with hyperbolic structure; on the other hand, the trefoil knot does not.
\end{remark}

Here, we present the statement of Mostow-Prasad rigidity in its most relevant form and case, with citation for source.

\begin{theorem}
\label{mostow-rigidity}
(Mostow's rigidity theorem for $3$-manifolds)  \cite[Theorem 8.3]{marden1974geometry}
Let $M, N$ be finite-volume hyperbolic $3$-manifolds. If $\pi_1(M)$ is isomorphic to $\pi_1(N)$, then $M$ is isometric to $N$.   
\end{theorem}

The following results and explanation on geometric topology of hyperbolic $3$-manifolds are adapted from Benedetti and Petronio's work in \cite{petronio1992lectures}, and have been altered to special cases; the results in the original text are more general.

\begin{definition}
Let $G$ be a group. We say that $G$ is \emph{virtually nilpotent} if $G$ contains some finite-index subgroup which is nilpotent.
\end{definition}

\begin{theorem}
(Margulis' Lemma) 
\cite[Adapted from Theorem D.1.1]{petronio1992lectures} There exists a constant $\eps$ such that for any properly discontinuous subgroup $\Gamma < Isom(\bb{H}^3)$ and any $x \in \bb{H}^3$, the group $\Gamma_\eps$ generated by the set $F_\eps(x) = \{\gamma \in \Gamma : d(x, \gamma(x)) \leq \eps\}$ is virtually nilpotent.
\end{theorem}

\begin{remark}
Adapted from \cite{petronio1992lectures} Let $M$ be a Riemannian manifold. If $\sigma$ is a piecewise differentiable path in $M$, we use $\ell(\sigma)$ to mean its length. For $\eps > 0$ we set $M_{(0, \eps]} = \{x \in M : \exists \langle \sigma \rangle \in \pi_1(M, x) \backslash \{1\}, \ell(\sigma) \leq \eps\}$ and $M_{[\eps, \infty]} = \{x \in M : \forall \langle \sigma \rangle \in \pi_1(M, x) \backslash \{1\}, \ell(\sigma) \geq \eps\}$.\\
Margulis' Lemma can be heuristically expressed in the following way: there exists a manifold-independent constant $\eps_3$ such that if $M$ is a complete oriented hyperbolic $3$-manifold and $x \in M$, the subgroup of $\pi_1(M, x)$ generated by $\eps_3$-short loops at $x$ is ``not very complicated''. This implies, more geometrically, that the $\eps_3$-thin part of $M$ is similarly ``not very complicated''.
\end{remark}

\begin{proposition}
\cite[Proposition D.3.12, paraphrasing]{petronio1992lectures}
All orientable finite-volume hyperbolic $3$-manifolds can be decomposed into a thick part, which is compact, and a thin part, which is composed of both cusps each of which is isotopic to $\bb{T}^2 \times \real^+$ with a distorted product metric, Euclidean in the first factor and exponentially decreasing in the second factor, and Margulis tubes, which are bounded neighborhoods of closed geodesics which are homeomorphic to the product of the geodesic with a $2$-ball. This is known as the \emph{thick-thin decomposition}.
\end{proposition}

Finally, we note that the groups $SL(2, \cpx)$ and $PSL(2, \cpx)$ show up frequently in this thesis. This is due to the fact that $Isom^+(\bb{H}^3) \cong PSL(2, \cpx)$: for an orientable finite volume hyperbolic $3$-manifold $M$, $\Gamma = \pi_1(M)$, we can define $M$ as $\bb{H}^3/\Gamma$, where $\Gamma$ is a subgroup of $Isom^+(\bb{H}^3) \cong PSL(2, \cpx)$, so that since $\widetilde{M} = \bb{H}^3$, $PSL(2, \cpx)$ itself acts by deck transformations on $M$. What is more, the inclusion representation $\rho: \Gamma \to PSL(2, \cpx)$ lifts to a representation $\widehat{\rho} : \Gamma \to SL(2, \cpx)$; for more on this, refer to MacLachan and Reid in \cite{aohm-book} on page 111-112.

\subsection{Representation Theory}
\label{repn-thy}

\begin{definition}
Let $\Gamma$ be a discrete subgroup of $SL(2, \cpx)$. The \emph{trace field} of $\Gamma$, $TF(\Gamma)$, is the field generated by all traces of all elements of $\Gamma$, which can also be written equivalently as $\q(tr\; \Gamma)$.
\end{definition}

We'll write $\q(tr\; \Gamma)$ in the few cases that we want to emphasize the trace field's nature as a number field, and $TF(\Gamma)$ otherwise.

\begin{theorem}
\cite[Theorem 3.1.2]{aohm-book} Let $M$ be a finite-volume orientable hyperbolic $3$-manifold, so that $\Gamma = \pi_1(M)$ is Kleinian. Then $\q(tr\; \Gamma)$ is a finite-degree extension.
\end{theorem}

\begin{definition}
The \emph{degree} of the trace field is the degree of the extension of $TF(\Gamma)$ over $\q$.
\end{definition}

\begin{proposition}
\label{D}
\cite[Proposition 1.5.2]{culler1983varieties} Let $\Pi$ be a finitely generated group, and let $\rho, \sigma$ be representations of $\Pi$ into $SL(2, \cpx)$ with corresponding characters $\chi_\rho, \chi_\sigma$. Let $\chi_\rho = \chi_\sigma$, and assume $\rho$ is an irreducible representation. Then $\rho, \sigma$ are conjugate.
\end{proposition}

\begin{remark}
We note that since traces (and thus characters) are invariant under conjugation, it suffices to consider the trace field of any choice of character within each conjugacy class. A natural way to think about the image of the character is therefore to look at its trace field, and the degree of the trace field over $\q$ is one of several natural measures of complexity.
\end{remark}

Next, we define notation to be used for the sake of clarity and concision for the rest of the thesis.

\begin{definition}
\label{charvar-notation}
Let $G$ be a group, $k$ a field, and $H$ a matrix group over $k$. Let $\mathrm{Hom_{Irr}}(G, H)$ be the set of irreducible representations of $G$ in $H$. We define
\begin{align*}
\chi^I_k(\Gamma) &:= \mathrm{Hom_{Irr}}(\Gamma, \mathrm{SL}(2, k))/\sim,
\end{align*}
where $\sim$ denotes the conjugacy relation. In particular,
\begin{align*}
\chi^I_\complex(\Gamma) &:= \mathrm{Hom_{Irr}}(\Gamma, \mathrm{SL}(2, \complex))/\sim,\\
\chi^I_{\finfieldcls{p}}(\Gamma) &:= \mathrm{Hom_{Irr}}(\Gamma, \mathrm{SL}(2, \finfieldcls{p}))/\sim.
\end{align*}
\end{definition}

This last we write as $\chi^I_p(\Gamma)$ for the sake of further concision. 

\begin{remark}
Theorem \ref{D} gives us a bijection between the set of irreducible characters, and the set of irreducible representations up to conjugacy. Accordingly, we call $\chi^I_\cpx(\Gamma)$ the $\cpx$-character variety of $\Gamma$. We note that character varieties truly are (reducible) varieties in their own right - for a detailed treatment, see \cite{culler1983varieties} and \cite{long2003integral}. 
\end{remark}

\subsection{The Sky Road and Contrapositive}
An important finiteness result for my thesis comes from Culler and Shalen. For clarity, I cite the result as stated in Reid and Wang \cite[Lemma 2.2]{reid1999non}. I have taken to calling the chain of results leading to it the ``Sky Road'', because when I first finally put the results together all in one place after poring over the half a dozen or so papers cited here, its trail from geometric topology through representation theory, algebraic geometry, and Bass-Serre theory, and back again to geometric topology, seemed almost magical to me. The pair of connected results at the end of the section are so vital to understanding my work that I have chosen to put a proof sketch of a weaker but more legible result here.\\
Let $M$ be an orientable hyperbolic $3$-manifold with $\Gamma = \pi_1(M)$, and suppose that $\Gamma$ has an affine algebraic curve $\mathscr{C}$ comprised of characters of representations $\rho: \Gamma \to SL(2, \cpx)$. Then we can complete that curve with ideal points of the variety by projectivizing. Next, given the ideal points of any variety, we can look at the localization of the rational function field to an ideal point, which will be a discrete valuation ring, and this gives us a valuation $v$ on $\Gamma$ itself. With a valuation $v$ on a group $\Gamma$, we can construct a tree $T_v$ which depends on $\Gamma, v$ in a natural way, and given that valuation tree, there exists a natural action $\Gamma \action T_v$ which lacks property FA, that is, it is free, lacks global fixed point, and without edge inversions. We can then take that action and quotient $T_v$ out by it, and get $T_v/\Gamma$ as a graph of groups. Pulling this back to $M$, the preimage of the edge midpoint gives us a locally-separating surface of $M$, which means that $M$ is Haken - for further details see Shalen \cite[Subsection 2.4]{shalen2002representations}.\\
There are a few subtleties here. First, we credit the path from curves of representations to valuation trees to the work of Culler and Shalen in \cite{culler1983varieties}, and the remaining part dealing with edge midpoint preimages and incompressible surfaces to Stallings in \cite{stallings1962fibering}. 

\begin{corollary}
\cite[Corollary 1.4.5]{culler1983varieties}
$\chi^I_\cpx(\Gamma)$ is a closed algebraic set. 
\end{corollary}
In \cite{culler1983varieties}, Culler and Shalen lay out the first half of the Sky Road in Section 1.3, moving from a curve of representations of a group $\Gamma = \pi_1(M)$, for $M$ a finite volume hyperbolic Haken 3-manifold, to the ideal point of the curve. Finally, they show in Theorem 2.2.1 that the localization of the rational function field on the curve to an ideal point yields a discrete valuation ring, which gives us a valuation on $\Gamma$ itself.
\begin{theorem}
\cite[Theorem 2.2.1]{culler1983varieties}
Let $C$ be a curve contained in $\chi_\cpx^I(\Gamma)$. Then to each ideal point $\widetilde{x}$ of $\widetilde{C}$ we can associate a nontrivial splitting of $\Gamma$.
\end{theorem}

On the other hand it is from Stallings's work in \cite{stallings1962fibering} and Shalen's additional explanation in \cite{shalen2002representations} that we get the rest of the Sky Road. Shalen spends Section 2 (by way of Bass-Serre theory) giving us a natural way to turn the valuation from the previous paper into a tree that $\Gamma$ acts on in a natural way freely, without global fixed point, and without edge inversions, and in particular in Subsection 2.4 he shows how preimages of edge midpoints yield pre-incompressible surfaces; that is, the preimage of any given edge midpoint is either incompressible or can be homotoped into incompressibility, since by the fact that we have a nontrivial action, at some point, whether through careful choice of midpoint or repeated homotopy and compression, we must finally find the bottom of the Haken hierarchy in a disjoint union of $3$-balls.\\
Fortunately for my work, we can take the contrapositive of every single one of these inferential steps. Let $N$ be an orientable non-Haken hyperbolic $3$-manifold with $\Delta = \pi_1(N)$. Since $N$ is non-Haken, it has no incompressible locally-separating surfaces, so that the preimage under the quotient map by $\Delta$ of any edge midpoint is always entirely compressible, and also that we can't generate a helpful graph of groups. Lacking that, we know that $\Delta$ has property FA, and in particular, the valuation tree $T_v$ doesn't exist, and thus neither does the valuation $v$. This means in turn that we have no ideal points to find, so that there exists no curve $\mathscr{C}$ of representations $\rho : \Delta \to SL(2, \cpx)$ to complete, so that the set of representations up to conjugacy must in fact be of dimension $0$. Since the set of representations is an algebraic variety of dimension $0$, it must in fact be finite.\\
Summarizing the above discussion, we deduce the following result:

\begin{theorem}
\label{A}
\cite[Lemma 2.2]{reid1999non}\\
Let $M$ be a non-Haken $3$-manifold of finite volume, with $\Gamma = \pi_1(M)$. Then $|\chi_\cpx^I(\Gamma)|$ is finite and $\Gamma$ has finitely many representations into $SL(2, \cpx)$ up to conjugacy.
\end{theorem}

\subsection{Profinite Groups}
We now introduce the other aspect of geometric group theory that we make use of in this thesis. For a more complete and formal treatment of the elementary characterizing properties of profinite completions of groups, I recommend reading through \cite{reid2013profinite}.

\begin{definition}
\label{inverse-sys}
Let $\{X_i\}_{i \in I}$ be a family of sets indexed by a partially-ordered set $(I, \leq)$, and suppose that for all $i \leq j, i, j \in I$, there is a function $f_{ij} : X_j \to X_i$ such that $f_{ii}$ is the identity, and for all $i \leq j \leq k$, $f_{jk} \circ f_{ij} = f_{ik}$. The \emph{inverse limit} of the $\{X_i\}$ is the subset of the direct product of the $\{X_i\}$ given by $X = \left\{\vec{x} \in \displaystyle\prod_{i \in I} X_i : \forall i, j \in I, i \leq j \Rightarrow \vec{x}_i = f_{ij}(\vec{x}_j)\right\}$, that is, the entries of $\vec{x}$ are compatible with all connection maps $f_{ij}$. If the $X_i$ are also topological spaces, the inverse limit inherits its topology from the product topology on $\prod X_i$, with the restriction that the projection maps must form commuting triangles as above.
\end{definition}

\begin{remark}
Any group can be made into a topological group by endowing it with the discrete topology.
\end{remark}

\begin{definition}
\label{profin-def}
A \emph{profinite} group is a topological group isomorphic to some inverse limit of an inverse system of finite groups with the discrete topology. Equivalently, profinite groups are compact, Hausdorff, totally disconnected topological groups, as per Ribes and Zalesskii in \cite[Theorem 1.1.12]{ribes2010profinite}.\\
The \emph{profinite completion} of a group $G$, which we denote by $\widehat{G}$, is that profinite group whose choice of finite groups is the set of all $G/N$, where $N$ ranges over the normal subgroups of $G$ of finite index, and the homomorphisms are given by the partial ordering of reverse containment of normal subgroups.
\end{definition}

We care about the profinite completion of a residually finite group $G$ primarily because it packages together all data on maps from $G$ to its quotient groups of finite order. 

\begin{lemma}
Let $G$ be a group, and let $\{N_i\}_{i \in I}$ range over all finite-index normal subgroups of $G$, where the indexing set $I$ is partially ordered by reverse set inclusion. Then $\{G/N_i\}_{i \in I}$ forms an inverse system.
\end{lemma}

\begin{proof}
We follow Definition \ref{inverse-sys} above. For $f_{ij} : G/N_j \to G/N_i$ we take the quotient map which looks, for $g \in G$, like $[g]_i \mapsto [g]_j$ on equivalence classes, which will make the $f_{ii}$ be the identity on $G/N_i$, and for $i \leq j \leq k \in I$, $f_{jk} \circ f_{ij} = f_{ik}$ as the composition of straightfoward quotient maps.
\end{proof}

\begin{definition}
\label{residual-finiteness}
Let $G$ be a group. We call $G$ \emph{residually finite} if for all $1 \neq g \in G$, there exists a finite group $H$ and a homomorphism $\phi : G \to H$ such that $\phi(g) \neq 1$.
\end{definition}

\begin{remark}
We may think of residual finiteness as the capacity for at least one of the finite-index (normal) subgroups of $G$ to tell an arbitrary $g \in G\backslash\{1\}$ apart from the identity, and we may think of the profinite completion of a group $\widehat{G}$ to be the packaging together of all of this finite-index information.
\end{remark}

\begin{lemma}
\cite[Adapted from Lemma 1.1.7]{ribes2010profinite}
Let $G$ be a group, with $\widehat{G}$ its profinite completion. Let $\iota : G \to \widehat{G}, g \mapsto ([g]_i)_{i \in I}$ be the canonical map sending $g$ to the $I$-indexed tuple of equivalence classes of $g$ under quotienting by the normal subgroups $\{N_i\}_{i \in I}$. Then $\iota$ has dense image.
\end{lemma}

The following proposition is well-known: see for example Reid \cite[Subsection 2.2]{reidicm} as well as Ribes and Zalesskii \cite[Page 78]{ribes2010profinite} .

\begin{proposition}
Given a group $G$, denote by $\mathcal{N}(G)$ the set of all finite-index normal subgroups of $G$. Then the following are equivalent:
\begin{itemize}
\item $G$ is residually finite;
\item $\displaystyle\bigcap_{N \in \mathcal{N}(G)} N = \{1\}$; and
\item The natural homomorphism $\iota : G \to \widehat{G}$ is injective.
\end{itemize}
\end{proposition}

\begin{lemma}
\cite[Lemma 3.2.1, p79]{ribes2010profinite} This canonical map $\iota$ satisfies the universal property that for any profinite group $H$ and group map $f : G \to H$, there exists a unique $g : \widehat{G} \to H$ such that $g \circ \iota = f$.
\end{lemma}

\begin{lemma}
Let $G$ be finitely generated. Since $SL(2, \overline{\finfield{p}})$ is locally finite, for any representation $\rho : G \to SL(2, \overline{\finfield{p}})$, $im\; \rho$ is also finite. 
\end{lemma}

\begin{corollary}
\label{local-finiteness}
Let $\rho : G \to SL(2, \overline{\finfield{p}})$ be a homomorphism. Then $\rho$ extends uniquely to a map $\widehat{\rho}: \widehat{G} \to SL(2, \overline{\finfield{p}})$.
\end{corollary}

\begin{definition}
\label{profin-equiv}
We say that two groups $G, H$ are \emph{profinitely equivalent} if $\widehat{G} \cong \widehat{H}$.
\end{definition}

\begin{definition}
\label{profin-rigid}
Let $G$ be a residually finite group. We say that a residually finite group $G$ is \emph{absolutely profinitely rigid} if for all residually finite groups $H$, whenever $\widehat{G} \cong \widehat{H}$, we have $G \cong H$. We say that a property is \emph{absolutely profinitely rigid} if any group possessing that property is profinitely rigid. We call two groups $G, G'$ \emph{absolutely profinitely distinguishable} if $\widehat{G} \not\cong \widehat{G'}$.
\end{definition}

\label{profin-rigid-2}
\begin{definition}
Let $G$ be a residually finite group, and $S$ a set of groups. We say that a residually finite group $G$ is \emph{relatively profinitely rigid} within $S$ if for all residually finite groups $H \in S$, whenever $\widehat{G} \cong \widehat{H}$, we have $G \cong H$. We say that a property is \emph{relatively profinitely rigid} within $S$ if any group in $S$ possessing that property is relatively profinitely rigid within $S$. We call two groups $G, G'$ \emph{relatively profinitely distinguishable} if $G, G' \in S$ and $\widehat{G} \not\cong \widehat{G'}$.
\end{definition}

\begin{remark}
The central focus of this paper is the profinite distinguishability of $3$-manifold groups. In the above definitions, we say that the rigidity or distinguishability is \emph{absolute} if we consider them in the general case, allowing $H$ to range over all residually finite groups, and we say that the rigidity or distinguishability is \emph{relative} to a specific family if we restrict our attention to a specific family of groups. In that case, we call the set of all groups profinitely equivalent to some group $G$ within that family the \emph{genus} of $G$ in that family.\footnote{The terminology here is not accidental, but borrowed from taxonomy, where family, genus, and species represent increasing levels of specificity of organism.} 
\end{remark}

\noindent In this paper, we restrict our attention to relative profinite rigidity within the category of $3$-manifold groups.

\section{Motivating Questions}

At this point, we can understand a few of the motivating open questions in geometric group theory. The questions here, which are foundational to the subfield, are adapted from Reid \cite[Section 4]{reidicm}.

\begin{question*}
Is $F_2$, the free group on two generators, absolutely profinitely rigid? In which families of groups is it relatively profinitely rigid?
\end{question*}

\begin{question*}
Let $\bb{S}$ be the set of $3$-manifold groups. Which $3$-manifold groups are relatively profinitely rigid within $\bb{S}$?
\end{question*}

\begin{question*}
Let $\bb{S}_n$ be the set of fundamental groups of hyperbolic $n$-manifolds. Which fundamental groups of hyperbolic $n$-manifolds are relatively profinitely rigid within $\bb{S}_n$?
\end{question*}

\begin{question*}
Let $M$ be a complete orientable finite-volume hyperbolic $3$-manifold. When is $\pi_1(M)$ absolutely profinitely rigid? When is $\pi_1(M)$ relatively profinitely rigid when $\bb{S}$ is the set of fundamental groups of hyperbolic $3$-manifolds?
\end{question*}

To answer these and other questions, new techniques will be required.

\section{Model Theory And Its Uses}
What could model theory be doing in a paper on geometric group theory and representation theory? Well, we use it as a means to prove the following theorem, whose proof can be found in Section \ref{more-model-theory}:	
\begin{restatable}{theorem}{passtofp}
\label{passtofp}
Suppose $|\chi^I_\complex(\Gamma)| = n$. Then for cofinitely many $p$, $|\chi^I_p(\Gamma)| = n$ as well.
\end{restatable}
Where model theory comes in is its ability to permit us to pass back and forth a given statement about $\cpx$ and the same single statement about $\finfieldcls{p}$ for cofinitely many primes $p$ - this is an example of a transfer principle. For early examples of this, see Lefschetz \cite{lefschetz-ag}, Weil \cite{weil-foundations}, or Seidenberg \cite{seidenberg1958comments}; for their context in a more general logical framework, see Ax \cite{ax1968elementary} in his proof of the Ax-Grothendieck Theorem, Barwise and Eklof \cite{barwise1969lefschetz}, or Robinson and Tarski \cite{robinson1953application}; and for the earliest modern form of the principle as used here, see Chapter 2 of Cherlin \cite{cherlin1976model}.\\
In any case, if for whatever reason you want to blackbox this, you can simply use the above Theorem \ref{passtofp}. Otherwise, however, we first need a few closely related definitions and a pair of results from model theory, which we'll detail later on in Section \ref{more-model-theory}; for a background reference not only for the model-theoretic techniques I use here but also introductory model theory as a whole, refer to \cite{marker1996introduction}.

\begin{definition}
The \emph{first-order theory of fields} has signature given by the constants $0, 1$, and the binary functions $+, \times$. It has the axioms that addition makes the set into an abelian group, multiplication is associative, commutative, and distributive with identity $1$, $\neg 0 = 1$, and removing $0$ makes the set into an abelian group under multiplication.
\end{definition}

\begin{definition}
The \emph{first-order theory of algebraically closed fields}, $ACF$, extends the first-order theory of fields by appending countably many axioms, one for each natural number, each of the form that every nontrivial polynomial of degree $n$ has at least one root.
\end{definition}

\begin{definition}
The \emph{first-order theory of algebraically closed fields of characteristic p}, $ACF_p$, extends $ACF$ by appending the additional axiom that $\overbrace{1 + 1 + \cdots + 1}^{p \mathrm{\;copies\;of\;} 1} = 0$.
\end{definition}

\begin{definition}
The \emph{first-order theory of algebraically closed fields of characteristic 0}, $ACF_0$, extends $ACF$ by appending countably many axioms, one for each prime $p$, each of the form that $\neg\overbrace{1 + 1 + \cdots + 1}^{p \mathrm{\;copies\;of\;} 1} = 0$.
\end{definition}

\section{Main Result}
\label{sec-main}

\subsection{More Model Theory}
\label{more-model-theory}
With the preliminaries well in hand, we can begin to discuss the specific way I apply tools from model theory to the study of profinite rigidity.

\begin{lemma}
\label{fpbar-and-c}
\cite[Corollary 1.2]{marker1996introduction}\\
Let $ACF_0$ be the first-order theory of algebraically closed fields of characteristic 0, and for rational prime $p$, let $ACF_p$ be the first-order theory of algebraically complete fields of characteristic $p$. Let $\Sigma$ be an $\scr{L}_1$-sentence. Then the following are equivalent:\\
(i) $ACF_0 \models \Sigma$;\\
(ii) $ACF_p \models \Sigma$ for cofinitely many choices of $p$;\\
(iii) $ACF_p \models \Sigma$ for infinitely many choices of $p$; and\\
(iv) $\cpx \models \Sigma$.
\end{lemma}


\begin{lemma}
\label{fields}
Let $S(k)$ be a first-order statement in the theory of a field $k$. Then $S(\finfieldcls{p})$ holds for infinitely (in fact, cofinitely) many choices of $p$ if and only if $S(\complex)$ holds.
\end{lemma}

\begin{proof}
Assume that $S(\cpx)$ holds. By Lemma \ref{fpbar-and-c}, this means that $ACF_0 \models S$, and thus also that $ACF_p \models S$ for cofinitely many $p$, that is, $S(\finfieldcls{p})$ holds for those $p$. On the other hand, assume that $S(\cpx)$ does not hold. Then by Lemma \ref{fpbar-and-c} again we must have $ACF_0 \models \neg S$, and thus also that $ACF_p \models \neg S$ for cofinitely many $p$, that is, $S(\finfieldcls{p})$ does not hold for those $p$.
\end{proof}

Whenever we want to leave the decision of which field we're working over until later, we will just write $ACF$. Fixing $k$ to be some arbitrary algebraically closed field\footnote{This still works for $k$ an arbitrary ring instead, but in that case, many of the properties below might be much weaker or have otherwise misleading names.}, we start by looking at how we can use $ACF$ to talk about matrices in $\mathrm{SL}(2, k)$. Let $x_1, x_2, x_3, x_4$ be variables in $ACF$. Then we define the predicate
\begin{equation*} M(x_1, x_2, x_3, x_4) \end{equation*} to be \begin{equation*} x_1 \cdot x_4 - x_2 \cdot x_3 = 1. \end{equation*}

The attentive reader may notice that this is exactly the defining relation for the determinant of a matrix in $\mathrm{M}(2, k)$ to be $1$ in terms of its elements. More subtly, and perhaps more powerfully, one may note a tactic that will be used throughout this section: namely, that we will make our predicates complex and full of equations, so that they can do the heavy lifting that a mere abstract tuple cannot do. More simply, though, we bundle $4$-tuples of variables that satisfy $M$ and notate them as matrices $A \in \mathrm{SL}(2, \complex)$ unless we really do need access to the entries.\\

To write that a given matrix is the identity is actually even easier. We define the predicate
\begin{equation*} Id(x_1, x_2, x_3, x_4) \end{equation*} to be \begin{equation*} x_1 = 1 \wedge x_2 = 0 \wedge x_3 = 0 \wedge x_4 = 1. \end{equation*}

We can also just write $Id(A)$, when we have a $4$-tuple as mentioned above. 

Before we can look at how to extend our method for talking about matrices of $\mathrm{SL}(2, k)$ in $ACF$ to a method for talking about representations $\Gamma \to \mathrm{SL}(2, k)$, we are going to need to be able to write predicates that verify that each relation of $\Gamma$ is satisfied. Consider how, for some finitely presented group \begin{equation*}\Gamma = \langle L | R \rangle,
\end{equation*} one might describe a new representation $\rho: \Gamma \to SL(2, k)$. It suffices to specify what the map does to a given choice of generators of $G$, $(l_i) \mapsto \rho(l_i)$. It is worth noting that this gives us a natural map $\mathrm{Hom}(\Gamma, SL(2, k)) \to k^{4l}$ determined by which element of $\mathrm{SL}(2, k)$ that particular $\rho \in \mathrm{Hom}(\Gamma, \mathrm{SL}(2, k))$ sends the ordered $l$-tuple of generators of $\Gamma$ to, interpreted by reading off the matrix entries. This is certainly injective - different elements of $\mathrm{Hom}(\Gamma, \mathrm{SL}(2, k))$ send at least one generator of $\Gamma$ to different matrices. What \cite{culler1983varieties} gives us is the conceptualization of the image as also some vanishing set $V_\Gamma(k)$, and then additionally that using the next few results (when we have proven them) that we can also can recognize which points of $k^{4l}$ are in $V_\Gamma(k)$ and are thus true representations by understanding $V_\Gamma(k)$ as the vanishing set of the polynomial relations in $R$, along with the polynomials ensuring that the generators map to elements of $\mathrm{SL}(2, k)$. But to make use of all thus, we have to tackle the challenge of how to communicate all of the machinery used here in $ACF$ first. Bringing this to $ACF$, let $(A_i) := A_1, \cdots, A_l$ be $4$-tuples such that

\begin{equation*} \displaystyle\bigwedge_{i=1}^{l} M(A_i). \end{equation*}\noindent That is, the vector in $k^{4l}$ can be thought of more helpfully as a $l$-tuple of $4$-tuples, each of which we conceptualize as a matrix. As such, we write the $l$-tuple $(A_i)_{i=1}^l$ as $\vec{A}$. A couple of lemmata on the relation between the entries of matrices and those of their products and inverses mean we can use $T$ to talk about the satisfaction of relations:

\begin{lemma}
\label{poly1}
Let $A, B \in \mathrm{SL}(2, k)$. Then the entries of $AB, A\inv$ are polynomial in the entries of $A, B$.
\end{lemma}


\begin{corollary}
\label{poly2}
Let $F \langle L \rangle$ be the set of freely reduced words on some finite set of letters and their inverses, $l = |L|$, $r \in F \langle L \rangle$. Let $f_r : k^{4l} \to k^4$ be the map treating successive $4$-tuples of the argument as matrix elements of a generating set, interpreting concatenation as matrix multiplication, and inverses of letters as inverses of generators, to take $r$ to its image under this choice of assignment. Then for all $\vec{z} \in k^{4l}$, $f_r(\vec{z})$ is polynomial in the $z_i$.
\end{corollary}


In particular, we use this in the case where the $z_i$ represent elements of $\mathrm{SL}(2, \complex)$.\\
Now that we have ensured that nothing goes horribly wrong when we talk about matrix inversions, matrix products, and free reduction in $\mathrm{SL}(2, k)$, we can now figure out how to use $ACF$ to write that a relationship $r \in R$ is satisfied by some $l$-tuple of matrices. For $r \in R$, we define the predicate

\begin{equation*} SAT_r(\vec{A}) \end{equation*} to be \begin{equation*} Id(f_r(\vec{A})). \end{equation*}\noindent That is, the $f_r$ map reads off the word $r$ and takes the $l$-tuple it spells out from the $A_i$, to the identity, as is required by the fact that the $r \in R$ are the relations of $G$.\\
Perhaps the most important immediate idea of this section is that we can use these previous predicates to write a statement in first-order logic that says whether or not a given $4l$-tuple is a representation of $G$ into $\mathrm{SL}(2, k)$.  In particular, we define \begin{equation*} REP_G(\vec{A}) \end{equation*} to be \begin{equation*} \left(\displaystyle\bigwedge_{i=1}^{l} M(A_i) \right) \wedge \left(\displaystyle\bigwedge_{r \in R} SAT_r(\vec{A}) \right), \end{equation*}
where the $G$ subscript reminds us that we started by fixing $G$ and a presentation for $G$ once and for all. Half-translating towards prose, this says that a $4l$-tuple $A$ corresponds to a representation of $G$ if all $l$ of the $4$-tuples are matrices of $\mathrm{SL}(2, k)$, and that all of the relations $r \in R$ of $G$ are satisfied. Later on we use this for discrete finite-covolume subgroups of $\mathrm{SL}(2, \cpx)$, designated by capital Greek letters such as $\Gamma$.

\begin{proposition}
\label{cs83}
\cite[Corollary to Proposition 1.4.1]{culler1983varieties}\\
Let $\Gamma$ be a finitely presented group, with finite presentation $\Gamma = \langle L | R \rangle$, and let $V_\Gamma(k) = \{\vec{A}\;|\;REP_\Gamma(\vec{A})\} = \{\vec{A}\;|\;V(f_r(\vec{A}) - I_2), M(A_i)_{i=1}^l\} \sbeq k^{4l}$. Then the map $\Phi_\Gamma : V_\Gamma \to \mathrm{Hom}(\Gamma, \mathrm{SL}(2, k))$ taking points in $V_\Gamma$ to the representations that they determine is well-defined and is a bijection.
\end{proposition}

\begin{remark}
The Hilbert Basis Theorem permits us to assume finite generation without loss of generality, rather than the stronger condition of finite presentation, and, in proofs, to pass to the case of finite presentation rather than finite generation.
\end{remark}

\begin{proof}
We start by noting that the first part of $REP_\Gamma$ uses the $M$ predicate to check that the $4l$-tuple genuinely is an $l$-tuple of matrices in $\mathrm{SL}(2, k)$. Thinking once more of $\vec{A}$ as a $4l$-tuple $v$, this defines a map $\psi_v: L \to \mathrm{SL}(2, k), a_i \mapsto A_i$. Then by the universal property of free groups, this extends to a map $\check{\psi}_v : F(L) \to \mathrm{SL}(2, k)$ under multiplication and inverses (see the diagram below).\\
Subsequently, the second part of $REP_\Gamma$ uses each $SAT(r)$ to check that each relation of $\Gamma$ is satisfied, and by the universal property of groups defined by presentations, this uniquely determines a representation $\rho_v: \Gamma \to \mathrm{SL}(2, k)$. But then that means that $\Phi_\Gamma(v) = \rho_v$.
\end{proof}

\[\begin{tikzcd}
	{L_\Gamma} && {F(L_\Gamma)} && \Gamma \\
	\\
	{\mathrm{SL}(2 , k)}
	\arrow["{\psi_v}", from=1-1, to=1-3]
	\arrow["{\check{\psi_v}}", from=1-3, to=1-5]
	\arrow[from=1-1, to=3-1]
	\arrow["{\rho_v}"{description}, dashed, from=1-5, to=3-1]
	\arrow[dashed, from=1-3, to=3-1]
\end{tikzcd}\]

\begin{remark}
After \cite{culler1983varieties}, we call $V_\Gamma(k)$ the \emph{$\mathrm{SL}(2, k)$-representation variety of $\Gamma$}.
\end{remark}

Now that we have established that we can talk about whether a given $4l$-tuple corresponds to a representation of $G$, we can talk about whether that representation is irreducible. As it turns out, though, it is much easier to start with reducibility. We recall that a representation into $\mathrm{SL}(2, k)$ is \emph{reducible} if the action by all of the images of the generators on $k^2$ fix some line through the origin: more formally, for some generator-dependent $\lambda \in k$,
\begin{equation*}
\left
(
\begin{array}{cc}
a & b \\
c & d
\end{array}
\right
)
\left
(
\begin{array}{c}
v_1 \\
v_2
\end{array}
\right
)
=
\left
(
\begin{array}{c}
\lambda v_1 \\
\lambda v_2
\end{array}
\right
),
\end{equation*}

where 
$\left
(
\begin{array}{cc}
a & b \\
c & d
\end{array}
\right
)$
ranges over the images of all generators of $\Gamma$. We can therefore define the predicate

\begin{equation*} RED(\vec{A}) \end{equation*} to be \begin{equation*} REP_G(\vec{A}) \wedge \exists a\; \exists b \displaystyle\bigwedge_{i = 1}^{n} \exists \lambda_i: ((a, b) \neq (0, 0)) \wedge A_i \cdot \langle a, b \rangle = \lambda_i \langle a, b \rangle, \end{equation*}
where we treat $\langle a, b \rangle$ as a column vector, and matrix and scalar multiplication are accordingly appropriately defined.

\begin{proposition}
For all $4l$-tuples $(A_i)$, $ACF \models RED(\vec{A})$ if and only if $\Phi(\vec{A})$ is a reducible representation.
\end{proposition}

We can then talk about irreducibility, defining the predicate 
\begin{equation*} IRREP(\vec{A}) \end{equation*} to be \begin{equation*} REP_G(\vec{A}) \wedge \lnot RED(\vec{A}). \end{equation*}

With a little more work, we can also talk about whether two representations are conjugate. We recall that two representations $\rho, \sigma : G \to \mathrm{SL}(2, k)$ are conjugate if there exists some matrix $M$ such that for all $i$, $M \rho(x_i) M\inv = \sigma(x_i)$, where $x_i$ is the $i^{th}$ generator of $G$ under some fixed choice of ordering. We can represent this in model theory by defining the predicate

\begin{equation*} CONJ(\vec{A}, \vec{B}) \end{equation*} to be \begin{equation*} REP_G(\vec{A}) \wedge REP_G(\vec{B}) \wedge \left(\exists C: M(C) \wedge \displaystyle\bigwedge_{j = 1}^{l} (C A_j C\inv = B_j)\right). \end{equation*}

\begin{proposition}
For all pairs of $4l$-tuples $\vec{A}, \vec{B}$, $ACF \models CONJ(\vec{A}, \vec{B})$ if and only if $\Phi(\vec{A}) \sim \Phi(\vec{B})$, that is, there exists some $C \in \mathrm{SL}(2, k)$ with $\tau_C \circ \Phi(\vec{A}) = \Phi(\vec{B})$, where $\tau_C$ is the inner automorphism of $\mathrm{SL}(2, k)$ that $C$ defines.
\end{proposition}

The whole point of this section, of course, was to be able to talk about the number of irreducible representations of a given finitely generated group $G$ into $\mathrm{SL}(2, k)$, up to conjugacy. However, this certainly is not a proper sentence in first-order logic. We might think of the ``half-translated'' version of $\Sigma_{G, n}$ as the following:

\begin{center}
There exist $n$ irreducible representations up to conjugacy of $G$ into $\mathrm{SL}(2, k)$, they are not conjugate to each other, and any other irreducible representation of $G$ into $\mathrm{SL}(2, k)$ must be conjugate to one of the $n$ representations previously mentioned.
\end{center}

We now have all the tools we need to write the previously mentioned sentence $\Sigma_{G, n}$; this will be almost exactly a predicate-by-predicate, symbol-by-symbol calque of what was written above as the half-translation, making use of the predicates we have constructed here. We write the first-order sentence $\Sigma_{G, n}$ as \begin{equation*} \displaystyle\bigwedge_{j = 1}^n \exists \vec{A}^{(j)}: IRREP(\vec{A}^{(j)}) \wedge \left( \displaystyle\bigwedge_{\substack{j, j' = 1\\j \neq j'}}^n \neg CONJ(\vec{A}^{(j)}, \vec{A}^{(j')})\right) \wedge \; \forall \vec{B}: \left(IRREP(\vec{B}) \Rightarrow \displaystyle\bigvee_{j=1}^n CONJ(\vec{B}, \vec{A}^{(j)})\right). \end{equation*}

\begin{theorem}
\label{sigma-gamma}
Let $G$ be a finitely generated group. Then there exists a sentence $\Sigma_{G, n}$ in $ACF$ such that for all algebraically closed fields $k$, $\Sigma_{G, n}(k)$ is true if and only if $|\chi_k(G)| = n$.
\end{theorem}

\begin{proof}
$(\Leftarrow)$ Suppose that $|\chi^I_\cpx(G)| = n$. Then by completeness of first-order logic, to verify that there exists a sentence $\Sigma_{G,n}$ in $ACF$ such that $ACF \models \Sigma_{G, n}$, it suffices to carefully step through the sentence itself to verify its meaning. We cut the sentence $\Sigma_{G, n}$ into three parts along the two conjunctions. The first clause asserts that there exists some family of $n$ tuples of appropriate length $\{\vec{A}^{(i)}\}_{j=1}^n$, each of which is an irreducible representation in itself. The second clause asserts that any distinct pair of those representations is nonconjugate. The final clause asserts that for all tuples $\vec{B}$ of the same length as the $\vec{A}^{(j)}$, if $B$ is an irreducible representation, then it must be conjugate to one of the $\vec{A}^{(j)}$. Taken in sum, the sentence asserts that $|\chi^I_\cpx(G)| = n$, and since we know it to be the case, $ACF$ models it.\\
$(\Rightarrow)$ Suppose that there exists a sentence $\Sigma_{G,n}$ in $ACF$ as constructed above such that $ACF \models \Sigma_{G, n}$. Since from the previous part we know that $\Sigma_{G, n}$ asserts that $|\chi^I_\cpx(G)| = n$ within $ACF$, by consistency of first-order logic we know that $|\chi^I_\cpx(G)| = n$.\\
Since $G$ is merely finitely generated and not finitely presented, we must invoke the Hilbert Basis Theorem: using it, we have that there exists a finitely presented $\tilde{G}, \eta : \tilde{G} \to G$ such that the induced map $\eta^* : \chi_\cpx^I(G) \to \chi_\cpx^I(\tilde{G})$ is a bijection. This $\tilde{G}$ is generated by any generating set for $G$, and its relations are the finitely many relations that correspond to the finitely many equations that the Hilbert Basis Theorem gives us. Then the sentence $\Sigma_{\tilde{G}, n}$ for $\tilde{G}$ also works for $G$.
\end{proof}

\begin{lemma}
\label{LA1}
Let $G$ be a finitely generated subgroup of $\mathrm{SL}(2, \finfieldcls{p})$. Then $|G|$ is finite.
\end{lemma}

\begin{proof}
To see this, we note that some generator will have an entry belonging to the largest $\finfield{p^k}$ among the set of generators, and neither addition nor matrix multiplication can increase the size; finally, any given $\mathrm{SL}(2, \finfield{p^k})$ is finite.
\end{proof}

\passtofp*

\begin{proof}
This follows immediately from Lemma \ref{fpbar-and-c} and Theorem \ref{sigma-gamma}; we recall that in particular, since $|\chi^I_\complex(\Gamma)| = n$, $|\chi^I_p(\Gamma)| = n$ for infinitely many $p$, and thus by Lemma \ref{fpbar-and-c}, for cofinitely many $p$.
\end{proof}

\begin{remark}
This will be the key tool allowing us to pass between representations into $\mathrm{SL}(2, \cpx)$, which we understand well, and $\mathrm{SL}(2, \finfieldcls{p})$, which we prize for its local finiteness. The goal will now be to use the theory of profinite distinguishability to prove the Main Theorem.
\end{remark}

\subsection{Constraints on Profinite Completions}
Now that we have established that we can pass (almost) freely between representations into $\mathrm{SL}(2, \cpx)$ and into $\mathrm{SL}(2, \finfieldcls{p})$, we can start to get a sense of what this means for the profinite distinguishability of groups.

\begin{lemma}
\label{profin-same-size-p}
Let $\Gamma, \Lambda$ be two finitely generated groups such that $\widehat{\Gamma} \cong \widehat{\Lambda}$. Suppose that $|\chi^I_p(\Gamma)| = n$ for cofinitely many $p$. Then $|\chi^I_p(\Gamma)| = |\chi^I_p(\Lambda)| = n$ for those $p$.
\end{lemma}

\begin{proof}
Consider the following commutative diagram.

\[\begin{tikzcd}
	{\Gamma} && {\mathrm{SL}(2, \overline{\bb{F}}_p)} \\
	\\
	{\widehat{\Gamma} = \widehat{\Lambda}} && {\Lambda}
	\arrow["{i_\Gamma}" description, from=1-1, to=3-1]
	\arrow["{i_\Lambda}" description, from=3-3, to=3-1]
	\arrow["{\widehat{\rho} = \widehat{\sigma}}" description, from=3-1, to=1-3]
	\arrow["{\rho}" description, from=1-1, to=1-3]
	\arrow["{\sigma}" description, from=3-3, to=1-3]
\end{tikzcd}\]

We note that by the universal property of profinite completions, any representation $\rho : \Gamma \to \mathrm{SL}(2, \finfieldcls{p})$ extends profinitely to a representation $\widehat{\rho} : \widehat{\Gamma} \to \mathrm{SL}(2, \finfieldcls{p})$ by the local finiteness of $\mathrm{SL}(2, \finfieldcls{p})$, as in Lemma \ref{LA1}. We may note that by the commutativity of the diagram, any representation from $\chi^I_p(\Gamma)$ must factor as a composition of the canonical injection into $\widehat{\Gamma}$ and a representation from $\chi^I_p(\widehat{\Gamma})$, so that $|\chi^I_p(\Gamma)| = |\chi^I_p(\widehat{\Gamma})|$. Looking at the other half of the diagram, we note that since $\widehat{\Gamma} = \widehat{\Lambda}$, the same argument applies in reverse: compositions of the canonical injection of $\Lambda$ into $\widehat{\Lambda}$ with representations from $\chi^I_p(\widehat{\Lambda})$ must yield all of $\chi^I_p(\Lambda)$, so that $|\chi^I_p(\Lambda)| = |\chi^I_p(\widehat{\Lambda})|$.
\end{proof}

\begin{theorem}
\label{profin-same-size-c}
Let $\Gamma, \Lambda$ be two finitely generated groups with $\widehat{\Gamma} \cong \widehat{\Lambda}$ and $|\chi^I_\cpx(\Gamma)| = n$ for some $n \in \nat$. Then  $|\chi^I_\cpx(\Lambda)| = |\chi^I_\cpx(\Gamma)| = n$.
\end{theorem}

\begin{proof}
Since $|\chi^I_\cpx(\Gamma)| = n$, by Lemma \ref{fields} and the statement of $\Sigma_{\Gamma, n}$, we know that $|\chi^I_p(\Gamma)| = n$ for cofinitely many $p$. By profinite equivalence and Lemma \ref{profin-same-size-p}, we know that $|\chi^I_p(\Lambda)| = |\chi^I_p(\Gamma)| = n$. Finally, by another application of Lemma \ref{fields} and the statement of $\Sigma_{\Lambda, n}$, $|\chi^I_\cpx(\Lambda)| = n$.
\end{proof}

\begin{theorem}
\label{repsgrow-c}
Let $\Gamma = \pi_1(M)$, where $M$ is a compact hyperbolic $3$-manifold. If $deg(TF(\Gamma)) \geq d$ for some $d \in \nat$, then $|\chi_\cpx^I(\Gamma)| \geq d$.
\end{theorem}

\begin{proof}
Let $\Gamma$ be the fundamental group of a finite-volume hyperbolic $3$-manifold, and let $k = TF(\Gamma)$ be its trace field. Let $\{\sigma_i\} : k \inj \complex$ be a family of distinct embeddings, and let $\theta : \Gamma \inj \mathrm{SL}(2, k)$ be the map realizing elements of $\Gamma$ as matrices. Let $\widehat{\sigma_i} : \mathrm{SL}(2, k) \inj \mathrm{SL}(2, \complex)$ be the map by interpretation of matrix entries induced by the $\{\sigma_i\}$. Then the result will follow if $\rho_i = \widehat{\sigma_i} \circ \theta$ is a representation $\rho_i : \Gamma \to \mathrm{SL}(2, \complex)$, and if for $i \neq j$, $\rho_i, \rho_j$ are nonconjugate.\\
To see this, recall that $k$ is a number field by \cite[Theorem 3.1.2]{aohm-book}. Denote its degree over $\q$ as $d = deg(k)$, so that $\{\sigma_i\}_{i=1}^d : k \inj \complex$ are its $d$-many embeddings. Then since the maps $\{\sigma_i\} : k \to \complex$ are all different, and are all embeddings, they cannot agree on every $\gamma \in \Gamma$: there must exist some $\gamma \in \Gamma$ such that $\widehat{\sigma_i} \circ \theta(\gamma) \neq \widehat{\sigma_j} \circ \theta(\gamma)$.\\
But then given that $\rho_i = \widehat{\sigma_i} \circ \theta$, $\rho_j = \widehat{\sigma_j} \circ \theta$, for $tr : \mathrm{SL}(2, \complex) \to \complex$ the trace map, $tr \circ \rho_i(\gamma) \neq tr \circ \rho_j(\gamma)$. Thus $\gamma$ represents a witnessing element of $\Gamma$ on which the representations $\rho_i, \rho_j$ have different traces, which by Proposition \ref{D} tells us that the two representations cannot be conjugate.
\end{proof}

Having shown that profinite equivalence means that the number of representations up to conjugacy into $\mathrm{SL}(2,k)$ (if finite) are the same between the profinitely equivalent groups, the goal is now to attack the main theorem.

\section{Main Theorem}
\label{main}
We need one last result from Long and Reid in \cite{long2003integral}.

\begin{theorem}
\label{B}
\cite[Theorem 3.2]{long2003integral}\\
Let $M$ be an orientable hyperbolic $3$-manifold of finite volume, and $d \in \nat$. Then there are finitely many pairs $(m, n)$ such that $tr(\rho(\pi_1(M_{m/n})) \in k$ where $deg(k/\q) \leq d$ as an extension.
\end{theorem}

\mainone*

\begin{proof}
It suffices to show that $\widehat{\Gamma} \cong \widehat{\Lambda}$ only for finitely many $\Lambda$; we thus assume that $\widehat{\Gamma} \cong \widehat{\Lambda}$. By assumption, $|\chi^I_\cpx(\Gamma)| < \infty$; let $|\chi^I_\cpx(\Gamma)| = d$. However, by Lemma \ref{profin-same-size-c}, we know that since $\widehat{\Gamma} \cong \widehat{\Lambda}$, $|\chi^I_\cpx(\Gamma)| = |\chi^I_\cpx(\Lambda)|$. Now, by Theorem \ref{B}, we know that there are at most finitely many choices of surgery coefficient resulting in a manifold with degree of trace field of fundamental group with at most a given degree $d + 1$, and by Theorem \ref{repsgrow-c}, we know that if $deg(TF(\Lambda)) > d$, then $|\chi^I_\cpx(\Lambda)| > d$ as well, so that it is exactly these finitely many choices of surgery coefficient where it is even possible for us to have $|\chi^I_\cpx(\Lambda)| = d$. Finally, by Proposition \ref{D}, we know that it suffices to check on the level of traces, since irreducible representations with the same trace are conjugate.
\end{proof}

The above result thus lends itself to the following more geometrically focused corollary:

\begin{corollary}
\label{cmainone}
Let $M$ be a one-cusped, finite-volume, hyperbolic 3-manifold. Suppose $M_{m/n}$ is a hyperbolic Dehn filling of $M$ with surgery coefficients $m/n$ and with finite character variety (for instance, a non-Haken such filling), and let $\Gamma = \pi_1(M_{m/n})$. Then for hyperbolic Dehn fillings with all but finitely many other choices of surgery coefficient $M_{m'/n'}$ with $\Lambda = \pi_1(M_{m'/n'})$, we get that $\widehat{\Gamma} \not\cong \widehat{\Lambda}$.
\end{corollary}

\begin{proof}
We may start by passing without loss of generality to the case where $M_{m'/n'}$ has hyperbolic structure, thanks to Theorems A and 8.4 from \cite{wilton2017distinguishing}. By assumption, $|\chi^I_\cpx(\Gamma)| < \infty$, so that Theorem \ref{mainone} applies.
\end{proof}

\begin{remark}
Liu in \cite{liu2020finitevolume} proves a more general version of this result using completely different methods.
\end{remark}

We can extend this to the question that actually motivated this entire line of inquiry with a little help from \cite{hatcher1982boundary}:

\begin{definition}
A knot $K$ is \emph{small} if its complement $S^3 \backslash K$ contains no closed incompressible surface.
\end{definition}

\begin{theorem}
\label{C}
\cite[Unnumbered Theorem from p. 373-374]{hatcher1982boundary}\\
Let $K$ be a small knot. Then all but finitely many of its Dehn fillings $M_{m/n} = \{S^3 \backslash K\}_{m/n}$ are non-Haken.
\end{theorem}

This allows us to narrow in on the trailhead for this line of thought: that we can use all of this to say something interesting about hyperbolic Dehn fillings of small knots.

\begin{corollary}
Let $K$ be a small knot such that $S^3 \backslash K = M$ is a one-cusped, finite-volume, hyperbolic 3-manifold. Let $\Gamma = \pi_1(M_{m/n})$ be the fundamental group of a non-Haken hyperbolic Dehn filling of $M$ with surgery coefficients $m/n$. Then for Dehn fillings with all but finitely many other choices of surgery coefficient $M_{m'/n'}$ with $\Lambda = \pi_1(M_{m'/n'})$, we get that $\widehat{\Gamma} \not\cong \widehat{\Lambda}$.
\end{corollary}

\begin{proof}
This follows from Theorems \ref{mainone} and \ref{C}, given the fact that knot complements have infinite cyclic first homology and thus that fillings whose slopes have different numerators have different first homology, which profinite completion detects.
\end{proof}

\section{Future Work}
\label{future}

My existing research has led me to a natural conjecture:
\begin{conjecture*}
In the main theorem above, we can drop the assumption that the character variety is finite, and consequently the implicit assumption that $M_{m/n}$ is non-Haken.
\end{conjecture*}

I propose to show this by working in $\overline{\q_p}$, figuring out a way to isolate zero-dimensional components of the character variety, showing that these are still isolated points, and showing that the number of these still increases without bound as the denominator of the surgery coefficient increases. I am currently in the middle of working on this, and I anticipate that further use of the kinds of model-theoretic tools I have already developed should help, as well as more classically algebrogeometric techniques. Further development along this line of inquiry would mean a deeper understanding of the nature of profinite flexibility of $3$-manifold groups, and possibly represent progress towards the discovery of a combined geometric/model-theoretic approach to open questions about absolute profinite rigidity. Although the work of Bridson et al \cite{bridson2018absolute} provides finitely many carefully-constructed examples of absolutely profinitely rigid groups, current approaches lack a means of generating infinite families of such groups. As-yet underdeveloped techniques grounded in model theory have shown promise in attacking that question.\\

In a similar vein, I have made a few conjectures and posed a few organizing questions regarding the future possibilities of my work.\\

\begin{conjecture*}
It is possible to recover the character variety $\chi^I_\cpx(\Gamma)$ from the profinite completion $\widehat{\Gamma}$.
\end{conjecture*}

\begin{conjecture*}
It is possible to detect and count the $0$-dimensional components of the character variety $\chi^I_\cpx(\Gamma)$ from the profinite completion $\widehat{\Gamma}$. 
\end{conjecture*}

\begin{question*}
Can the positive-dimensional components of the character variety also be detected and counted all together? Can they then be separated by dimension?
\end{question*}

\begin{conjecture*}
It is possible, through the use of model theoretic techniques, to construct infinite families of absolutely profinitely rigid groups.
\end{conjecture*}

Notably for this last conjecture, my work thus far has only used the most basic of model-theoretic tricks and my own extensions of those initial tricks, and has been the first to combine these two perspectives; I plan to investigate whether further and more powerful model-theoretic tools might be waiting to be applied to geometric group theory.\\

Nonetheless, even very simple questions of absolute profinite ridigity remain. For instance, given some arbitrary $3$-manifold group, it is in general an open question as to whether it is absolutely profinitely rigid. Additionally, as in Noskov, Remeslennikov, and Roman'kov in \cite{noskov1979infinite}, it is a long-standing and foundational open question whether there exists any residually finite group profinitely equivalent to the free group on two generators $F_2$.

\bibliographystyle{abbrv}
\bibliography{References}

\begin{thebibliography}{10}

\bibitem{agol2013virtual}
I.~Agol.
\newblock The virtual {H}aken conjecture.
\newblock {\em Doc. Math}, 18:1045--1087, 2013.
\newblock With an appendix by Agol, Daniel Groves, and Jason Manning.

\bibitem{ax1968elementary}
J.~Ax.
\newblock The elementary theory of finite fields.
\newblock {\em Annals of Mathematics}, pages 239--271, 1968.

\bibitem{barwise1969lefschetz}
J.~Barwise and P.~Eklof.
\newblock Lefschetz's principle.
\newblock {\em Journal of Algebra}, 13(4):554--570, 1969.

\bibitem{bridson2018absolute}
M.~R. Bridson, D.~McReynolds, A.~W. Reid, and R.~Spitler.
\newblock Absolute profinite rigidity and hyperbolic geometry.
\newblock {\em Annals of Mathematics}, 192:679--719, 2020.

\bibitem{cherlin1976model}
G.~L. Cherlin.
\newblock Model theoretic algebra.
\newblock {\em The Journal of Symbolic Logic}, 41(2):537--545, 1976.

\bibitem{culler1983varieties}
M.~Culler and P.~B. Shalen.
\newblock Varieties of group representations and splittings of 3-manifolds.
\newblock {\em Annals of Mathematics}, pages 109--146, 1983.

\bibitem{gromov1981hyperbolic}
M.~Gromov.
\newblock Hyperbolic manifolds according to {T}hurston and {J}{\o}rgensen.
\newblock In {\em S{\'e}minaire Bourbaki vol. 1979/80 Expos{\'e}s 543--560},
  pages 40--53. Springer, 1981.

\bibitem{hatcher1982boundary}
A.~Hatcher.
\newblock On the boundary curves of incompressible surfaces.
\newblock {\em Pacific Journal of Mathematics}, 99(2):373--377, 1982.

\bibitem{hatcher-notes}
A.~Hatcher.
\newblock Notes on basic 3-manifold topology.
\newblock \url{https://pi.math.cornell.edu/~hatcher/3M/3Mdownloads.html}, 2007.

\bibitem{hempel2004three}
J.~Hempel.
\newblock {\em 3-Manifolds}, volume 349.
\newblock American Mathematical Soc., 2004.

\bibitem{lefschetz-ag}
S.~Lefschetz.
\newblock {\em Algebraic geometry}.
\newblock Princeton University Press, 2015.

\bibitem{liu2020finitevolume}
Y.~Liu.
\newblock Finite-volume hyperbolic 3-manifolds are almost determined by their
  finite quotient groups.
\newblock arXiv preprint arXiv: 2011.09412, 2020.

\bibitem{long2003integral}
D.~Long and A.~W. Reid.
\newblock Integral points on character varieties.
\newblock {\em Mathematische Annalen}, 325(2):299--321, 2003.

\bibitem{aohm-book}
C.~Maclachlan and A.~W. Reid.
\newblock {\em The Arithmetic of Hyperbolic 3-Manifolds}.
\newblock Graduate Texts in Mathematics. Springer, 2003.

\bibitem{maher2010random}
J.~Maher.
\newblock Random {H}eegaard splittings.
\newblock {\em Journal of Topology}, 3(4):997–1025, 2010.

\bibitem{malcev1940isomorphic}
A.~Malcev.
\newblock On isomorphic matrix representations of infinite groups.
\newblock {\em Matematicheskii Sbornik}, 50(3):405--422, 1940.

\bibitem{marden1974geometry}
A.~Marden.
\newblock The geometry of finitely generated kleinian groups.
\newblock {\em Annals of Mathematics}, 99(3):383--462, 1974.

\bibitem{marker1996introduction}
D.~Marker, M.~Messmer, and A.~Pillay.
\newblock {\em Introduction to the model theory of fields}, volume~5 of {\em
  Lecture {N}otes in {L}ogic}.
\newblock Springer-Verlag, Berlin, 1996.

\bibitem{milnor1962unique}
J.~Milnor.
\newblock A unique decomposition theorem for 3-manifolds.
\newblock {\em American Journal of Mathematics}, 84(1):1--7, 1962.

\bibitem{noskov1979infinite}
G.~A. Noskov, V.~N. Remeslennikov, and V.~A. Roman'kov.
\newblock Infinite groups.
\newblock {\em Itogi Nauki i Tekhniki. Seriy ``Algebra. Geometriya.
  Topologiya"}, 17:65--157, 1979.

\bibitem{petronio1992lectures}
R.~Petronio, R.~Benedetti, and C.~Petronio.
\newblock {\em Lectures on Hyperbolic Geometry}.
\newblock Universitext (Berlin. Print). Springer Berlin Heidelberg, 1992.

\bibitem{przytycki2018mixed}
P.~Przytycki and D.~Wise.
\newblock Mixed 3-manifolds are virtually special.
\newblock {\em Journal of the American Mathematical Society}, 31(2):319--347,
  2018.

\bibitem{reid2013profinite}
A.~W. Reid.
\newblock Profinite properties of discrete groups.
\newblock In {\em Proceedings of Groups St Andrews}, pages 73--104, 2013.

\bibitem{reidicm}
A.~W. Reid.
\newblock Profinite rigidity.
\newblock {\em Proceedings of the International Congress of Mathematics (Rio de
  Janeiro)}, 2:1211--1234, 2018.

\bibitem{reid1999non}
A.~W. Reid and S.~Wang.
\newblock Non-{H}aken 3-manifolds are not large with respect to mappings of
  non-zero degree.
\newblock {\em Communications in Analysis and Geometry}, 7(1):105--132, 1999.

\bibitem{ribes2010profinite}
L.~Ribes and P.~Zalesskii.
\newblock {\em Profinite Groups}.
\newblock Ergebnisse der Mathematik und ihrer Grenzgebiete. 3. Folge / A Series
  of Modern Surveys in Mathematics. Springer Berlin Heidelberg, 2010.

\bibitem{robinson1953application}
A.~Robinson and A.~Tarski.
\newblock On the application of symbolic logic to algebra.
\newblock {\em Journal of Symbolic Logic}, 18(2), 1953.

\bibitem{seidenberg1958comments}
A.~Seidenberg.
\newblock Comments on lefschetz's principle.
\newblock {\em The American Mathematical Monthly}, 65(9):685--690, 1958.

\bibitem{shalen2002representations}
P.~B. Shalen et~al.
\newblock Representations of 3-manifold groups.
\newblock {\em Handbook of geometric topology}, pages 955--1044, 2002.

\bibitem{stallings1962fibering}
J.~Stallings.
\newblock On fibering certain 3-manifolds, topology of 3-manifolds and related
  topics (proc. the univ. of georgia institute, 1961).
\newblock {\em Prentice-Hall, Englewood Cliffs, NJ VA Steklov Mathematical
  Institute of Russian Academy of Science}, 8:95--100, 1962.

\bibitem{thurston1982three}
W.~P. Thurston.
\newblock Three dimensional manifolds, {K}leinian groups and hyperbolic
  geometry.
\newblock {\em Bulletin of the American Mathematical Society}, 6(3):357--381,
  1982.

\bibitem{weil-foundations}
A.~Weil.
\newblock {\em Foundations of algebraic geometry}, volume~29.
\newblock American Mathematical Society, 1946.

\bibitem{wilton2017distinguishing}
H.~Wilton and P.~Zalesskii.
\newblock Distinguishing geometries using finite quotients.
\newblock {\em Geometry \& Topology}, 21(1):345--384, 2017.

\bibitem{wilton2019profinite}
H.~Wilton and P.~Zalesskii.
\newblock Profinite detection of 3-manifold decompositions.
\newblock {\em Compositio Mathematica}, 155(2):246--259, 2019.

\bibitem{wise2009research}
D.~T. Wise.
\newblock Research announcement: the structure of groups with a quasiconvex
  hierarchy.
\newblock {\em Electron. Res. Announc. Math. Sci}, 16(44-55):6, 2009.

\end{thebibliography}
\end{document}